\documentclass[11pt]{article}
\usepackage[margin=1in]{geometry} 
\usepackage{amssymb,amsfonts,amsmath,bbm,mathrsfs,stmaryrd}
\usepackage{xcolor}
\usepackage{url}

\usepackage{enumerate}

\usepackage[colorlinks,
             linkcolor=black!75!red,
             citecolor=blue,
             pdftitle={},
             pdfproducer={pdfLaTeX},
             pdfpagemode=None,
             bookmarksopen=true
             bookmarksnumbered=true]{hyperref}

\usepackage{tikz}
\usetikzlibrary{arrows,calc,decorations.pathreplacing,decorations.markings,intersections,shapes.geometric,through,fit,shapes.symbols,positioning,decorations.pathmorphing}

\usepackage{braket}

\usepackage[amsmath,thmmarks,hyperref]{ntheorem}
\usepackage{cleveref}

\creflabelformat{enumi}{#2(#1)#3}

\crefname{section}{Section}{Sections}
\crefformat{section}{#2Section~#1#3} 
\Crefformat{section}{#2Section~#1#3} 

\crefname{subsection}{\S}{\S\S}
\AtBeginDocument{%
  \crefformat{subsection}{#2\S#1#3}%
  \Crefformat{subsection}{#2\S#1#3}%
}

\crefname{subsubsection}{\S}{\S\S}
\AtBeginDocument{%
  \crefformat{subsubsection}{#2\S#1#3}%
  \Crefformat{subsubsection}{#2\S#1#3}%
}

\theoremstyle{plain}

\newtheorem{lemma}{Lemma}[section]

\newtheorem{corollary}[lemma]{Corollary}
\newtheorem{theorem}[lemma]{Theorem}

\theoremstyle{nonumberplain}

\theoremstyle{plain}
\theorembodyfont{\upshape}
\theoremsymbol{\ensuremath{\blacklozenge}}

\newtheorem{definition}[lemma]{Definition}
\newtheorem{example}[lemma]{Example}
\newtheorem{remark}[lemma]{Remark}

\crefname{definition}{definition}{definitions}
\crefformat{definition}{#2definition~#1#3} 
\Crefformat{definition}{#2Definition~#1#3} 

\crefname{ex}{example}{examples}
\crefformat{example}{#2example~#1#3} 
\Crefformat{example}{#2Example~#1#3} 

\crefname{remark}{remark}{remarks}
\crefformat{remark}{#2remark~#1#3} 
\Crefformat{remark}{#2Remark~#1#3} 

\crefname{convention}{convention}{conventions}
\crefformat{convention}{#2convention~#1#3} 
\Crefformat{convention}{#2Convention~#1#3} 

\crefname{notation}{notation}{notations}
\crefformat{notation}{#2notation~#1#3} 
\Crefformat{notation}{#2Notation~#1#3} 

\crefname{table}{table}{tables}
\crefformat{table}{#2table~#1#3} 
\Crefformat{table}{#2Table~#1#3}

\crefname{lemma}{lemma}{lemmas}
\crefformat{lemma}{#2lemma~#1#3} 
\Crefformat{lemma}{#2Lemma~#1#3} 

\crefname{proposition}{proposition}{propositions}
\crefformat{proposition}{#2proposition~#1#3} 
\Crefformat{proposition}{#2Proposition~#1#3} 

\crefname{corollary}{corollary}{corollaries}
\crefformat{corollary}{#2corollary~#1#3} 
\Crefformat{corollary}{#2Corollary~#1#3} 

\crefname{theorem}{theorem}{theorems}
\crefformat{theorem}{#2theorem~#1#3} 
\Crefformat{theorem}{#2Theorem~#1#3} 

\crefname{enumi}{}{}
\crefformat{enumi}{(#2#1#3)}
\Crefformat{enumi}{(#2#1#3)}

\crefname{assumption}{assumption}{Assumptions}
\crefformat{assumption}{#2assumption~#1#3} 
\Crefformat{assumption}{#2Assumption~#1#3} 

\crefname{equation}{}{}
\crefformat{equation}{(#2#1#3)} 
\Crefformat{equation}{(#2#1#3)}

\numberwithin{equation}{section}

\theoremstyle{nonumberplain}
\theoremsymbol{\ensuremath{\blacksquare}}

\newcommand\pf[1]{\newtheorem{#1}{Proof of \Cref{#1}}}

\newcommand\bK{{\mathbb K}}

\newcommand\cA{{\mathcal A}}
\newcommand\cB{{\mathcal B}}
\newcommand\cC{{\mathcal C}}

\newcommand\cE{{\mathcal E}}

\newcommand\cM{{\mathcal M}}

\newcommand\cO{{\mathcal O}}

\newcommand\cR{{\mathcal R}}
\newcommand\cS{{\mathcal S}}
\newcommand\cT{{\mathcal T}}

\newcommand\cX{{\mathcal X}}
\newcommand\cY{{\mathcal Y}}

\DeclareMathOperator{\id}{id}

\newcommand{\cat}[1]{\textsc{#1}}

\newcommand\spr[1]{\cite[\href{https://stacks.math.columbia.edu/tag/#1}{Tag {#1}}]{stacks-project}}
\newcommand{\qedhere}{\mbox{}\hfill\ensuremath{\blacksquare}}

\title{Centers of categorified endomorphism rings}
\author{Alexandru Chirvasitu}

\begin{document}

\date{}

\newcommand{\Addresses}{{
  \bigskip
  \footnotesize

  \textsc{Department of Mathematics, University at Buffalo, Buffalo,
    NY 14260-2900, USA}\par\nopagebreak \textit{E-mail address}:
  \texttt{chirvasitua@gmail.com}

}}

\maketitle

\begin{abstract}
  We prove that for a large class of well-behaved cocomplete categories $\mathcal C$ the weak and strong Drinfeld centers of the monoidal category $\mathcal{E}$ of cocontinuous endofunctors of $\mathcal{C}$ coincide. This generalizes similar results in the literature, where $\mathcal{C}$ is the category of modules over a ring $A$ and hence $\mathcal{E}$ is the category of $A$-bimodules.
\end{abstract}


\vspace{.5cm}



\section*{Introduction}

The present note is motivated by the following result from \cite{acg} (see Theorem 2.10 therein):

\begin{theorem}\label{th:init}
  For a ring $A$, the weak and strong centers of the monoidal category ${}_A\cM_A$ of $A$-bimodules coincide.
\end{theorem}

We give a refresher on the terminology in \Cref{subse:cent} below, pausing here only for a broad-strokes perspective on the result.

As seen from \Cref{def:cent,def:wcent} below, \Cref{th:init} says, essentially, that a certain morphism
\begin{equation*}
  A\otimes V\to V\otimes A
\end{equation*}
of $A\otimes A$-bimodules (for a bimodule $V\in {}_A\cM_A$ underlying a weak-center object) is automatically an isomorphism. The proofs of \cite[Propositions 2.5 and 2.6]{acg} make it clear that this is the type of rigidity phenomenon familiar from the theory of {\it descent} in ring theory and / or algebraic geometry \cite{ko}. In the latter setup one typically starts with commutative rings $R\to S$ and an $S$-module $M$ and seeks to recover an $R$-module $M_R$ such that
\begin{equation*}
  M\cong S\otimes_R M_R;
\end{equation*}
in other words, the goal is to {\it descend} the $S$-module to an $R$-module. The sort of structure necessary to achieve this in good cases (e.g. $S$ is faithfully flat over $R$ \cite[Th\'eor\`eme 3.2]{ko}) is a {\it descent datum} (see \cite[discussion preceding Proposition 3.1]{ko}): an $S\otimes S$-module morphism
\begin{equation*}
  g:S\otimes M\to M\otimes S
\end{equation*}
(where `$\otimes$' means `$\otimes_R$') such that
\begin{enumerate}[(1)]
\item the diagram
  \begin{equation*}
    \begin{tikzpicture}[auto,baseline=(current  bounding  box.center)]
      \path[anchor=base] 
      (0,0) node (l) {$S\otimes S\otimes M$}
      +(3,.5) node (u) {$S\otimes M\otimes S$}
      +(6,0) node (r) {$M\otimes S\otimes S$}
      ;
      \draw[->] (l) to[bend left=6] node[pos=.5,auto] {$\scriptstyle g_{23}$} (u);
      \draw[->] (u) to[bend left=6] node[pos=.5,auto] {$\scriptstyle g_{12}$} (r);
      \draw[->] (l) to[bend right=6] node[pos=.5,auto,swap] {$\scriptstyle g_{13}$} (r);
    \end{tikzpicture}
  \end{equation*}
  commutes (with the indices indicating the tensorands on which $g$ operates), and
\item the morphism
  \begin{equation*}
    M\to S\otimes M\stackrel{g}{\longrightarrow} M\otimes S\to M
  \end{equation*}
  is the identity, where the leftmost arrow is the natural inclusion obtained by tensoring the unit $R\to S$ of $S$ with $\id_M$ and the rightmost arrow is multiplication by scalars in $S$.
\end{enumerate}

Under these circumstances it turns out \cite[Proposition 3.1]{ko} that in fact $g$ is automatically an isomorphism. This is essentially the same phenomenon as that captured in \Cref{th:init} in the broader context of non-commutative rings. 

In attempting to isolate precisely what it is about categories of bimodules that occasions such rigidity results one is led to consider the celebrated Eilenberg-Watts theorem (\cite[Theorem 1]{wts} or \cite{eil}):
\begin{itemize}
\item ${}_A\cM_A$ is equivalent to the category of cocontinuous (i.e. colimit-preserving \cite[\S V.4]{macl}) endofunctors of the category ${}_A\cM$ of left $A$-modules (or its right-handed version $\cM_A$),
\item such that the monoidal structure given by `$\otimes_A$' is identified with endofunctor composition.
\end{itemize}

This is the starting point for the generalization of \Cref{th:init} appearing as \Cref{th:main} below. The pattern we extrapolate can be summarized as follows (with a forward reference to \Cref{subse:2alg} below for category-theoretic terminology).

\begin{itemize}
\item One can substitute other ``well-behaved'' cocomplete categories $\cC$ for ${}_A\cM$;
\item and their {\it duals} $\cC^*\cong $ consisting of cocontinuous functors $\cC\to\text{(some ``base'' category)} $ for $\cM_A$;
\item and their {\it endomorphism 2-rings}
  \begin{equation}\label{eq:initbx}
    \cE:=\cC\boxtimes \cC^*\cong \text{cocontinuous endofunctors of }\cC
  \end{equation}
  for ${}_A\cM_A$. 
\end{itemize}

For \Cref{eq:initbx} to be both meaningful and valid $\cC$ needs to be what in \Cref{def:dlz} (and \cite[Definition 1.1]{us}) we refer to as {\it dualizable} (this is what `well-behaved' means in the above discussion). With all of this behind us, \Cref{th:main} reads more or less as follows.

\begin{theorem}
  If $\cC$ is a dualizable locally presentable category then the weak center of its category of cocontinuous endofunctors coincides with its strong center.
\end{theorem}

In addition to recovering \Cref{th:init}, this applies to categories $\cC$ going beyond modules, as we recall in \Cref{se:main}: $\cC$ can be, for instance,
\begin{itemize}
\item the category $\cM^C$ of right-comodules over a {\it right-semiperfect} (\cite[p.369]{MR0498663}) coalgebra over a field;
\item the category $\cat{Qcoh}([X/G])$ of quasicoherent sheaves over the quotient stack $[X/G]$ where $X$ is affine and $G$ is a {\it virtually linearly reductive} \cite[\S 1]{MR1395068} linear algebraic group acting on $X$.
\end{itemize}

\subsection*{Acknowledgements}

This work was partially supported by NSF grant DMS-2001128. 

I am grateful for illuminating comments by Ana Agore and for the referee's careful reading and insightful comments.

\section{Preliminaries}\label{se.prel}

Some standard background on monoidal categories is needed, as covered for instance in \cite[Chapter XI]{kas}, \cite[\S 5.1]{cp-q}, \cite[Chapter XI]{macl}, and any number of other sources.

\subsection{Some 2-algebra}\label{subse:2alg}

We reprise some terminology from \cite[\S 2]{cf}.

\begin{definition}
  \begin{enumerate}[(a)]
  \item A {\it 2-abelian group} is a locally presentable category in the sense of \cite[Definition 1.17]{ar}. 2-abelian groups form a 2-category \cat{2Ab} with left adjoints as 1-morphisms and natural transformations as 2-morphisms. 
  \item A {\it 2-ring} is a 2-abelian group $\cC$ which is in addition a monoidal category with tensor product `$\otimes$', so that all functors of the form $x\otimes -$ and $-\otimes x$ are left adjoints. 2-rings similarly form a 2-category \cat{2Rng} with {\it monoidal} left adjoints as 1-morphisms. 
  \item A {\it commutative 2-ring} is a 2-ring additionally equipped with a symmetry (i.e. it is a symmetric monoidal category). As before, these form the 2-category \cat{2ComRng} with symmetric monoidal left adjoints as 1-morphisms.
  \end{enumerate}
\end{definition}

It turns out (e.g. \cite[Corollary 2.2.5]{cf}) that \cat{2Ab} is symmetric monoidal, being equipped with a tensor product denoted by `$\boxtimes$'. For 2-abelian groups $\cA$ and $\cB$ their tensor product $\cA\boxtimes \cB$ is the universal recipient of a bifunctor
\begin{equation*}
  \cA\times \cB\to \cA\boxtimes \cB
\end{equation*}
that is separately cocontinuous (i.e. a ``bilinear map'' of 2-abelian groups). The symmetric monoidal structure lifts to \cat{2Rng} and \cat{2ComRng} in the sense that if $\cA$ and $\cB$ are 2-rings so is $\cA\boxtimes \cB$ in a natural fashion, etc.

This machinery allows us to employ the usual language of rings and modules in the context of 2-abelian groups:

\begin{definition}
  Let $\cR$ be a 2-ring. A {\it left (2-)$\cR$-module} is a 2-abelian group $\cX$ equipped with a morphism $\cR\boxtimes \cX\to \cX$ in \cat{2Ab}, satisfying the obvious unitality and associativity conditions. {\it Right} (2-)$\cR$-modules are defined analogously, as are bimodules, etc.

  The respective 2-categories of left or right or bimodules are denoted by ${}_{\cR}\cM$, $\cM_{\cR}$ and ${}_{\cR}\cM_{\cS}$. respectively.
\end{definition}

As usual, we have tensor product 2-bifunctors
\begin{equation*}
  {}_{\cR}\cM_{\cS}\times   {}_{\cS}\cM_{\cT} \stackrel{\boxtimes_{\cS}}{\longrightarrow}   {}_{\cR}\cM_{\cT}:
\end{equation*}
For 2-$\cS$-modules $\cX$ and $\cY$ with module-structure functors
\begin{equation*}
  \cX\boxtimes \cS\stackrel{\triangleleft}{\longrightarrow}\cX
  \quad
  \text{and}
  \quad
  \cS\boxtimes \cY\stackrel{\triangleright}{\longrightarrow}\cY
\end{equation*}
we can define $\cX\boxtimes_{\cS}\cY$ as the universal functor $\pi$ admitting a natural isomorphism $\theta$ as depicted below:
\begin{equation*}
  \begin{tikzpicture}[auto,baseline=(current  bounding  box.center)]
    \path[anchor=base] 
    (0,0) node (l) {$\cX\boxtimes \cS\boxtimes \cY$}
    +(3,.5) node (u) {$\cX\boxtimes \cY$}
    +(3,-.5) node (d) {$\cX\boxtimes \cY$}
    +(6,0) node (r) {$\cX\boxtimes_{\cS} \cY$}
    +(3,0) node (c) {$\scriptstyle \cong\ \Downarrow \theta$}
    ;
    \draw[->] (l) to[bend left=6] node[pos=.5,auto] {$\scriptstyle \triangleleft\boxtimes\id$} (u);
    \draw[->] (u) to[bend left=6] node[pos=.5,auto] {$\scriptstyle \pi$} (r);
    \draw[->] (l) to[bend right=6] node[pos=.5,auto,swap] {$\scriptstyle \id\boxtimes\triangleright $} (d);
    \draw[->] (d) to[bend right=6] node[pos=.5,auto,swap] {$\scriptstyle \pi$} (r);
  \end{tikzpicture}
\end{equation*}
This is a 2-colimit, and can be obtained as a combination of a {\it co-inserter} and a {\it co-equifier}: constructions dual to the inserters and equifiers of \cite[\S 2.71 and Lemma 2.76]{ar}, for instance. \cite[Proposition 2.1.11]{cf} outlines a more explicit construction for such 2-colimits in \cat{2Ab}, which do exist.

In particular, for a {\it commutative} 2-ring $\cR$, the 2-category ${}_{\cR}\cM\cong \cM_{\cR}$ is symmetric monoidal under $\boxtimes_{\cR}$.

\begin{definition}\label{def:2ralg}
  For a commutative 2-ring $\cR$ an {$\cR$-algebra} (or 2-$\cR$-algebra for extra precision) is an algebra in the symmetric monoidal 2-category ${}_{\cR}\cM$. 
\end{definition}

It turns out that \cat{2Ab} is not only symmetric monoidal but also {\it monoidal-closed}, i.e. has {\it internal homs}. More precisely, we have the familiar hom-tensor adjunction in the present higher-categorical setting (see e.g. \cite[\S 6.5]{kel-enr}, \cite[Exercise 1.l]{ar}). The following result is a ``relative'' version of \cite[Lemma 2.7]{us} (which cites the preceding two sources), in the sense that it deals with modules over 2-rings rather than plain 2-abelian groups. The techniques involved in the proofs are no different.

\begin{lemma}
  Let $\cR$, $\cS$ and $\cT$ be three 2-rings.
  \begin{enumerate}[(a)]
  \item For any two bimodules
    \begin{itemize}
    \item $\cX\in {}_{\cR}\cM_{\cT}$
    \item $\cY\in {}_{\cS}\cM_{\cT}$
    \end{itemize}
    the category
    \begin{equation*}
      \cat{Hom}_{\cT}(\cY,\cX):=\{\text{left adjoints }\cX\to \cY\text{ compatible with the 2-module structures}\}
    \end{equation*}
    has a natural structure of a $\cR$-$\cS$-bimodule. 
  \item This gives us, for each bimodule $\cY\in {}_{\cS}\cM_{\cT}$, a 2-adjunction
    \begin{equation*}
 \begin{tikzpicture}[auto,baseline=(current  bounding  box.center)]
   \path[anchor=base] 
   (0,0) node (1) {${}_{\cR}\cM_{\cS}$} 
   +(6,0) node (2) {${}_{\cR}\cM_{\cT}$}
   ;
   \draw[->] (1) to[bend left=6] node[pos=.5,auto] {$\scriptstyle -\boxtimes_{\cS}\cY$} (2);
   \draw[->] (2) to[bend left=6] node[pos=.5,auto] {$\scriptstyle \cat{Hom}_{\cT}(\cY,-)$} (1);
 \end{tikzpicture}
\end{equation*}
with the top arrow as the left (2-)adjoint. 
\end{enumerate}
  \qedhere
\end{lemma}

\begin{remark}
  We leave it to the reader to formulate analogous versions for tensoring on the left rather than right, etc.
\end{remark}

\begin{definition}\label{def:dlz}
  Let $\cR$ be a 2-ring and $\cX$ a left $\cR$-module.
  \begin{enumerate}[(a)]
  \item The {\it dual} $\cX^*$ of $\cX$ over $\cR$ is $_{\cR}\cat{Hom}(\cX,\cR)$; it is a {\it right} $\cR$-module. Similarly, duals of right modules are naturally left modules.
  \item If $\cR$ is commutative the 2-$\cR$-module $\cX$ is {\it (1-)dualizable over $\cR$} if the canonical morphism
    \begin{equation}\label{eq:can}
      \cX\boxtimes_{\cR}\cX^*\stackrel{\cat{can}}{\longrightarrow} \cat{End}_{\cR}(\cX)
    \end{equation}
    is an isomorphism of 2-$\cR$-modules. 
  \end{enumerate}
\end{definition}

\begin{remark}
  For any 2-ring $\cR$ and 2-$\cR$-module $\cX$ $\cat{End}_{\cR}(\cX)$ is naturally a 2-ring (and a $2$-$\cR$-algebra when $\cR$ is commutative), with composition as the tensor product and $\id_{\cX}$ as the unit.
\end{remark}

Dualizable objects (typically over $\cR=\cat{Vect}_{\bK}$ for some field $\bK$) were the focus of \cite{us}, where we give alternative characterizations of dualizability in \cite[Lemma 3.1]{us}. In particular, it is enough to require that the identity
\begin{equation*}
  \id_{\cX}\in \cat{End}_{\cR}(\cX)
\end{equation*}
belong to the image of \Cref{eq:can}.

\subsection{Centers}\label{subse:cent}

Recall (e.g. \cite[Definition XIII.4.1]{kas} or \cite[Definition 3]{js-tort}):

\begin{definition}\label{def:cent}
  Let $(\cC,\otimes,{\bf 1})$ be a monoidal category. The {\it (Drinfeld) center} $Z(\cC)$ of $\cC$ is the category of pairs $(x,\theta)$ where $x\in \cC$ is an object and
  \begin{equation}\label{eq:brd}
    \theta: -\otimes x\stackrel{\cong}{\longrightarrow}x\otimes -
  \end{equation}
  is a natural isomorphism satisfying the following conditions (suppressing the associativity constraints in the monoidal category):
  \begin{enumerate}[(1)]
  \item\label{item:1} For $y,z\in \cC$ the diagram
    \begin{equation*}
      \begin{tikzpicture}[auto,baseline=(current  bounding  box.center)]
        \path[anchor=base] 
        (0,0) node (l) {$y\otimes z\otimes x$}
        +(3,.5) node (u) {$y\otimes x\otimes z$}
        +(6,0) node (r) {$x\otimes y\otimes z$}
        ;
        \draw[->] (l) to[bend left=6] node[pos=.5,auto] {$\scriptstyle \id_y\otimes \theta_z$} (u);
        \draw[->] (u) to[bend left=6] node[pos=.5,auto] {$\scriptstyle \theta_y\otimes\id_z$} (r);
        \draw[->] (l) to[bend right=6] node[pos=.5,auto,swap] {$\scriptstyle \theta_{y\otimes z}$} (r);
      \end{tikzpicture}
    \end{equation*}
    commutes, and
  \item\label{item:2} the isomorphism
    \begin{equation*}
      \theta_{\bf 1}:{\bf 1}\otimes x\to x\otimes {\bf 1}
    \end{equation*}
    is the canonical one attached to the monoidal structure $(\cC,\otimes,{\bf 1})$. 
  \end{enumerate}
\end{definition}

\begin{remark}
  In fact, in \Cref{def:cent} condition \Cref{item:2} follows from \Cref{item:1}, but this uses the fact that $\theta$ is an {\it isomorphism}; we have displayed both conditions with an eye towards \Cref{def:wcent} below.
\end{remark}

Following \cite[Definition 4.3]{sch-dbl} (where the notion seems to have been introduced) and \cite[\S 1.1]{acg}, we give

\begin{definition}\label{def:wcent}
  For $(\cC,\otimes,{\bf 1})$ as in \Cref{def:cent}  the {\it weak right center} $WZ_r(\cC)$ is the category of pairs $(x,\theta)$ as above, satisfying conditions \Cref{item:1} and \Cref{item:2}, but requiring only that \Cref{eq:brd} be a natural transformation.

  One defines the weak {\it left} center $WZ_{\ell}(\cC)$ analogously, requiring a natural transformation
  \begin{equation*}
    x\otimes -\to -\otimes x
  \end{equation*}
  instead.

  Unless specified otherwise {\it weak center} means weak {\it right} center, and we simply write $WZ$ for $WZ_r$.
\end{definition}

\section{Main results}\label{se:main}

\begin{theorem}\label{th:main}
  Let $\cR$ be a commutative 2-ring, $\cX\in {}_{\cR}\cM$ a dualizable $\cR$-module, and
  \begin{equation}\label{eq:endrng}
    \cE:=\cX\boxtimes_{\cR} \cX^* \cong \cat{End}_{\cR}(\cX)
  \end{equation}
  its endomorphism ring. Then, the canonical fully faithful inclusion
  \begin{equation}\label{eq:z-to-wz}
    Z(\cE)\to WZ(\cE)
  \end{equation}
  is an equivalence.
\end{theorem}

This requires some preliminary discussion and tooling, starting with the observation that this weak-equals-strong principle cannot be expected to hold in general: \Cref{def:wcent} is indeed a weakening of \Cref{def:cent}, in the sense that there are examples of monoidal categories (even $2$-algebras) $\cC$ where the canonical functor \Cref{eq:z-to-wz} is not an equivalence.

\begin{example}\label{ex:simple}
  Let $(\cC,\otimes,{\bf 1})$ be any symmetric, additive monoidal category. Recall (e.g. \cite[Definition 2.1]{helm-cent}) that its {\it Bernstein center} is the (commutative) ring of natural endomorphisms of the identity functor. Any element $\theta$ of the Bernstein center gives an element $({\bf 1}, \psi)$ of the weak center $WZ(\cC)$, whereby
  \begin{equation*}
    y\cong y\otimes{\bf 1}\stackrel{\psi_y}{\longrightarrow} {\bf 1}\otimes y\cong y
  \end{equation*}
  is simply $\theta_y$, provided $\theta_{\bf 1}$ is the identity. Furthermore, if $\theta_y$ fails to be an isomorphism for any $y$ we obtain an element outside the plain center $Z(\cC)$.

  All of this is easily arranged. Let $\cC$, for instance, be the (symmetric, monoidal) category $\mathrm{Rep}_f(G)$ of finite-dimensional complex representations over a finite group $G$. Every $y\in \mathrm{Rep}_f(G)$ has a canonical ${\bf 1}$-{\it isotypic component}: the space $y^G$ of all $G$-invariant vectors in $y$. There is an element $\theta$ of the Bernstein center that surjects every object onto this isotypic component:
  \begin{equation*}
    y\ni v\stackrel{\theta_y}{\longmapsto} \frac 1{|G|}\sum_{g\in G}gv\in y
  \end{equation*}
  If $G$ is non-trivial then $\theta_y$ will fail to be an isomorphism on those $y$ that are not sums of copies of ${\bf 1}$, and we have an example as required above.
\end{example}

Since $\cR$ is our ``base ring'' throughout the discussion we henceforth abbreviate `$\boxtimes_{\cR}$' to simply `$\boxtimes$', and similarly for $\cat{Hom}:=\cat{Hom}_{\cR}$. Recall also our notation \Cref{eq:endrng} for the endomorphism 2-ring $\cE$ of $\cX$.

Because $\cX$ is assumed dualizable over $\cR$, the canonical morphism
\begin{equation*}
  \cX\to \cX^{**}
\end{equation*}
is an isomorphism (of abelian 2-groups, i.e. an equivalence of categories). It follows from this that $\cE$ is also dualizable and in fact self-dual, and we can identify
\begin{equation}\label{eq:allsame}
  \cat{End}_{\cR}(\cE)\ \cong\
  \cE\boxtimes \cE^*\ \cong\
  \cE\boxtimes \cE\ \cong\
  \cX\boxtimes \cX^* \boxtimes \cX\boxtimes \cX^*. 
\end{equation}
Given that $\cat{2Ab}$ is a {\it symmetric} monoidal 2-category, there is some choice in how we identify the right and left-hand sides of \Cref{eq:allsame}. In the sequel, it will be convenient to make this identification by pairing the two {\it middle} tensorands on the right-hand side of \Cref{eq:allsame} against $\cE\cong \cX\boxtimes \cX^*$ in the obvious fashion (by pairing each $\cX$ to a $\cX^*$). Concretely, simple-tensor object
\begin{equation*}
  x\boxtimes f\boxtimes y\boxtimes g \in \cX\boxtimes \cX^* \boxtimes \cX\boxtimes \cX^*
\end{equation*}
corresponds to the element
\begin{equation*}
  \cE\ni \psi\longmapsto f(\psi(y)) x\boxtimes g\in \cX\boxtimes \cX^*\cong \cE
\end{equation*}
of $\cat{End}_{\cR}(\cE)$.

Now fix an object $(e,\theta)\in WZ(\cE)$ of the weak right center and consider the two endomorphisms
\begin{equation*}
  -\otimes e \text{ and } e\otimes -\quad \in\quad \cat{End}(\cE) \cong \cE\boxtimes \cE.
\end{equation*}
With the above convention in mind, they are identifiable, respectively, with
\begin{equation}\label{eq:outid}
  {\bf 1}\boxtimes e\quad\text{and}\quad e\boxtimes {\bf 1}
\end{equation}
where
\begin{equation*}
  {\bf 1}:=\id_{\cX}={\bf 1}_{\cE}\in \cE = \cat{End}(\cX)\cong \cX\boxtimes \cX^*
\end{equation*}
is the identity functor on $\cX$ (i.e. the monoidal unit of $\cE$). We caution the reader that the tensor product in \Cref{eq:outid} is {\it external}, i.e. it is not to be confused with the internal tensor product `$\otimes$' of $\cE$. Indeed, under the latter we of course have
\begin{equation*}
  {\bf 1}\otimes e\cong e\cong e\otimes{\bf 1}
\end{equation*}
(as in any monoidal category).

The natural transformation
\begin{equation*}
 -\otimes e \to e\otimes -
\end{equation*}
that constitutes the structure of a weak-center element (\Cref{def:wcent}) translates to a morphism 
\begin{equation}\label{eq:newtheta}
  {\bf 1}\boxtimes e \stackrel{\theta}{\longrightarrow} e\boxtimes{\bf 1}
\end{equation}
in $\cE\boxtimes \cE$ (denoted slightly abusively by the same symbol `$\theta$' we used for the natural transformation in \Cref{def:wcent}). The conditions \Cref{item:1,item:2} can then be recast in terms of \Cref{eq:newtheta} as we explain presently.

To express condition \Cref{item:1} we need to work in the {\it triple} tensor product $\cE^{\boxtimes 3}$. To that end, we consider morphisms between tensor products of $e$ and two copies of ${\bf 1}$, with two indices among $1$, $2$ and $3$ indicating where $\theta$ operates. Thus:
\begin{align*}
  \theta_{12}:=\theta\boxtimes \id_{\bf 1}&:{\bf 1}\boxtimes e\boxtimes{\bf 1} \to e\boxtimes{\bf 1}\boxtimes{\bf 1},\\
  \theta_{23}:=\id_{\bf 1}\boxtimes \theta &:{\bf 1}\boxtimes {\bf 1}\boxtimes e \to {\bf 1}\boxtimes e\boxtimes{\bf 1},\\
\end{align*}
and similarly,
\begin{equation*}
  \theta_{13}:
  {\bf 1}\boxtimes {\bf 1}\boxtimes e
  \to
  e\boxtimes{\bf 1}\boxtimes{\bf 1}
\end{equation*}
is the morphism acting identically on the middle tensorand and as $\theta$ on the two outer ones. \Cref{item:1} in \Cref{def:cent} can now be recovered simply as
\begin{equation}\label{eq:thetasdesc}
  \theta_{13} = \theta_{12}\circ \theta_{23}:
  {\bf 1}\boxtimes {\bf 1}\boxtimes e
  \to
  e\boxtimes{\bf 1}\boxtimes{\bf 1}.  
\end{equation}

Next, denote by
\begin{equation*}
  m:\cE\boxtimes \cE \to \cE
\end{equation*}
the ``multiplication'' morphism, imparting on $\cE=\cat{End}(\cX)$ its monoidal category structure. In terms of the decomposition
\begin{equation*}
  \cE\boxtimes \cE\cong \cX\boxtimes \cX^*\boxtimes \cX\boxtimes \cX^*
\end{equation*}
$m$ is simply the evaluation of the two middle tensorands $\cX$ and $\cX^*$ against each other. With this in place, condition \Cref{item:2} in \Cref{def:cent} simply asks that
\begin{equation*}
  m(\theta):{\bf 1}\otimes e\to e\otimes{\bf 1} 
\end{equation*}
be the canonical isomorphism, i.e. the identity once we have made the usual identifications 
\begin{equation*}
  {\bf 1}\otimes e\cong
  e\cong
  e\otimes{\bf 1}.  
\end{equation*}
In short, for future reference:
\begin{equation}\label{eq:mtheta}
  m(\theta) = \id_{e}:e\cong {\bf 1}\otimes e\to e\otimes{\bf 1}\cong e. 
\end{equation}

\pf{th:main}
\begin{th:main}
  Since we already know that \Cref{eq:z-to-wz} is fully faithful (as is immediate from \Cref{def:cent,def:wcent}), it remains to show that it is essentially surjective: for an arbitrary object $(e,\theta)\in WZ(\cE)$ the morphism \Cref{eq:newtheta} is an isomorphism in $\cE\boxtimes \cE$. What we will in fact do is identify the inverse of $\theta$: it is precisely
  \begin{equation*}
    \theta':=\tau\circ\theta\circ\tau:e\boxtimes{\bf 1}\to {\bf 1}\boxtimes e,
  \end{equation*}
  where $\tau$ is the tensorand-reversal functor on $\cE\boxtimes \cE$.

  Denote by
  \begin{equation*}
    m_{13}:\cE\boxtimes \cE\boxtimes \cE\to \cE\boxtimes \cE
  \end{equation*}
  the functor that multiplies the outer (first and third) tensorands of the domain onto the second tensorand of the codomain. We then have
  \begin{align*}
    \theta\ &=\ m_{13}(\theta_{23})\quad \text{and}\\
    \theta'\ &=\ m_{13}(\theta_{12}),
  \end{align*}
  meaning that 
  \begin{equation}\label{eq:1313}
    m_{13}(\theta_{13}) = m_{13}(\theta_{12}\circ \theta_{23}) = \theta'\circ \theta
  \end{equation}
  (where the first equality uses \Cref{eq:thetasdesc} above). On the other hand though, \Cref{eq:mtheta} implies that the left-hand side $m_{13}(\theta_{13})$ of \Cref{eq:1313} is nothing but the identity, and thus
  \begin{equation*}
    \theta'\circ\theta = \id_{{\bf 1}\boxtimes e}. 
  \end{equation*}
  The other composition $\theta\circ\theta'$ is treated similarly, so we do not repeat the argument.
\end{th:main}

\begin{remark}\label{re:steal}
  The proof of \Cref{th:main} given above is a paraphrase, in the present categorified context, of an argument familiar from descent theory. See e.g. \cite[Proposition 3.1]{ko}. Where $\cX$ would have been the category of modules over (in those authors' notation) a ring $S$.
\end{remark}

In the context of $k$-linear 2-abelian groups (for some field $k$) the examples of dualizable $\cC$ in \cite{us} are all, abstractly, of the form
\begin{equation}\label{eq:pshv}
  \text{$k$-linear functors }\Gamma^{op}\to {}_k\cat{Vect}
\end{equation}
for small $k$-linear categories $\Gamma$. These are also
\begin{itemize}
\item the $k$-linear abelian categories admitting a generating set of {\it small} projective objects \cite[\S 3.6, Corollary 6.4]{pop};
\item the $k$-linear locally presentable admitting a strongly generating set of small projective objects \cite[Theorem 5.26]{kel-enr};
\end{itemize}
recall that the small projective objects in an abelian category $\cC$ are simply those $x\in \cC$ for which the representable functor
\begin{equation*}
  \mathrm{hom}(x,-):\cC\to \cat{Set}
\end{equation*}
is cocontinuous. This is taken as the definition of the hyphenated term `{\it strong-projective}' in \cite[\S 5.5]{kel-enr}, and we reuse that term in \Cref{cor:smproj} below for consistency.

Conversely, \cite[Lemma 3.5]{us} shows that {\it all} categories of the form \Cref{eq:pshv} are dualizable 2-modules over ${}_k\cat{Vect}$. The same argument goes through for arbitrary commutative 2-rings $\cR$ (in place of ${}_k\cat{Vect}$), so we have

\begin{corollary}\label{cor:smproj}
  Let $\cR$ be a commutative 2-ring and $\cX$ a 2-$\cR$-module with a strong generating set of small-projective objects. Then, the weak center of the monoidal category
  \begin{equation*}
    \cX\boxtimes_{\cR}\cX^*\cong \cat{End}_{\cR}(\cX)
  \end{equation*}
  coincides with its strong center.
  \qedhere
\end{corollary}

We end with some examples of categories falling under the scope of \Cref{cor:smproj} (and hence \Cref{th:main}).

\begin{example}\label{ex:coalg}
Throughout the present discussion we assume $C$ is a coalgebra over a field. 
  
By \cite[Theorem 1.3]{us}, \Cref{th:main} applies to categories of right comodules $\cM^C$ over {\it right-semiperfect} coalgebras $C$ in the sense of \cite[p.369]{MR0498663}: every right $C$-comodule has a projective cover.

This of course includes cosemisimple coalgebras (i.e. those with only projective modules or equivalently, direct sums of simple coalgebras; \cite[Definition, p.290]{MR0252485} or \cite[Definition 2.4.1]{mont}). 
\end{example}

\begin{example}\label{ex:stk}
  Overlapping \Cref{ex:coalg} to a degree, consider a linear algebraic group \cite[\S 1.6]{bor} $G$ acting on an affine scheme $X$ and the category
  \begin{equation}\label{eq:quotst}
    \cat{Qcoh}(X)^G\cong \cat{Qcoh}([X/G])
  \end{equation}
  of $G$-equivariant quasicoherent sheaves on $X$, or equivalently, as \Cref{eq:quotst} recalls (\spr{06WV}), that of quasicoherent sheaves on the quotient stack \spr{044O} $[X/G]$.

  According to \cite[Theorem 1.5]{us} \Cref{eq:quotst} is dualizable provided $G$ is {\it virtually linearly reductive} in the sense of \cite[\S 1]{MR1395068}: $G$ has a normal linearly reductive closed algebraic subgroup $H\trianglelefteq G$ such that $G/H$ is a finite group scheme. This means that
  \begin{itemize}
  \item the Hopf algebra $\cO(H)$ is cosemisimple while $\cO(G/H)$ is finite-dimensional;
  \item equivalently by \cite[Theorem, p.76]{MR1395068}, the Hopf algebra $\cO(G)$ of regular functions on $G$ is (left and right) semiperfect in the sense of \Cref{ex:coalg}.
  \end{itemize}  
\end{example}

\begin{example}\label{ex:genalg}
  One can generalize \Cref{ex:stk} as follows. Note that the category \Cref{eq:quotst} can be recovered as that of $\cO(X)$-modules (where $\cO(X)$ is the algebra of regular functions on $X$) internal to the category of $\cO(G)$-comodules. In short:
  \begin{equation*}
    \cat{Qcoh}([X/G])\cong \cM_{\cO(X)}^{\cO(G)}. 
  \end{equation*}
  Mimicking this construction, we can take $\cX$ in \Cref{th:main} to be the category $\cM_A$ of (right, say) modules over an algebra $A$ internal to the commutative 2-algebra $\cR$ (so that in \Cref{ex:stk} we would have $\cR=\cM^{\cO(G)}$ and $A=\cO(X)$).

  This means that \Cref{th:main} applies, for instance, to categories of graded modules over graded algebras (for arbitrary grading monoids), etc.
\end{example}



\addcontentsline{toc}{section}{References}

\Addresses

\end{document}